# A MAX-AUTOREGRESSIVE MODEL WITH MAX-SEMI-STABLE MARGINALS

S. Satheesh[*1] and E. Sandhya[**]

[*]NEELOLPALAM, S. N. Park Road, Trichur – 680 004, India.

*ssatheesh1963@yahoo.co.in*

[**]Department of Statistics, Prajyoti Niketan College, Pudukkad, Trichur – 680 301, India.

*esandhya@hotmail.com*

**Abstract**

*The structure of stationary first order max-autoregressive schemes with max-semi-stable marginals is studied. A connection between semi-selfsimilar extremal processes and this max-autoregressive scheme is discussed resulting in their characterizations. Corresponding cases of max-stable and selfsimilar extremal processes are also discussed.*

## 1. Introduction

This investigation is motivated by Satheesh and Sandhya (2007) wherein the structure of the first order autoregressive (AR(1)) model with semi-stable marginals were studied. Max-AR(1) models were introduced and characterized in Satheesh and Sandhya (2006) using max-semi-selfdecomposable (SD) laws. Here we need the following notions.

**Definition.1.1** (Megyesi, 2002). A non-degenerate *d.f* $F$ is max-semi-stable($a,b$) if

$$F(x) = exp\{-x^{-\alpha}h(ln(x))\},\ x>0,\ \alpha>0, \qquad (1.1)$$

where $h(x)$ is a positive bounded periodic function with period $ln(b)$, $b>1$, and there exists an $a>1$ such that $ab^{-\alpha} = 1$. This class the extended Frechet, is denoted by $\Phi_{\alpha,a,b}$, or

$$F(x) = exp\{-|x|^{\alpha}h(ln(|x|))\},\ x<0,\ \alpha>0, \qquad (1.2)$$

where $h(x)$ is a positive bounded periodic function with period $|ln(b)|$, $b<1$, and there exists an $a>1$ such that $ab^{\alpha} = 1$. This class the extended Weibull, is denoted by $\Psi_{\alpha,a,b}$. When $h(x)$ is a constant (1.1) and (1.2) represent max-stable *d.f*s.

**Definition.1.2** (Becker-Kern, 2001). A non-degenerate *d.f* $F$ is max-semi-SD($c$) if for some $c>1$ and $v \in \mathbf{R}$ there is a non-degenerate *d.f* $H$ such that

$$F(x) = F(c^v x + \beta) H(x) \; \forall \; x \in \mathbf{R}, \tag{1.3}$$

where $\beta=0$ *if* $v \neq 0$ and $\beta = ln(c)$ *if* $v=0$. If (1.3) holds for every $c>1$, then $F$ is max-SD.

---

[1] *Author for correspondence.*

Extremal Processes (EP) $\{Y(t), t \geq 0\}$ are processes with increasing right continuous sample paths and independent max-increments. In extreme value theory EPs $\{Y(t)\}$, with stationary max-increments and $Y(o) = 0$ form the counterpart of Levy processes. The univariate marginals of an EP determine its finite dimensional distributions. We will refer to EPs whose univariate marginals are max-semi-stable as max-semi-stable($a,b$) EPs. Selfsimilar (SS) processes are those that are invariant in distribution under suitable scaling of time and space. Pancheva (1998, 2000) has developed SS (semi-SS) EPs and showed that EPs whose univariate marginals are max-stable (max-semi-stable) constitute a class of SS (semi-SS) EPs having stationary max-increments. In general EPs are discussed with state space in $\mathbf{R}^d$, $d>1$ integer since their max-increments are max-infinitely divisible (max-ID) and all *d.f*s on $\mathbf{R}$ are max-ID. However, we can discuss max-stable, max-semi-stable, max-SD and max-semi-SD laws on $\mathbf{R}$ and so we restrict our discussion to *d.f*s on $\mathbf{R}$ and the invariance *w.r.t* linear normalization. Thus on $\mathbf{R}$ we have;

**Definition.1.3** An EP $\{Y(t)\}$ is semi-SS *if* for some $b>0$ there is an exponent $H \in \mathbf{R}$ such that

$$\{Y(bt)\} \stackrel{d}{=} \{b^H Y(t)\}. \tag{1.4}$$

In terms of the *d.f* $G(u)$ of $Y(1)$ this is equivalent to $\{G(b^H u)\}^{bt} = \{G(u)\}^t$. If (1.4) holds for any $b>0$ then $\{Y(t)\}$ is SS. We write $\{Y(t)\}$ is $(b,H)$-semi-SS ($\{Y(t)\}$ is $H$-SS).

From Satheesh and Sandhya (2006) we have the following definition and results.

**Definition.1.4** A sequence $\{X_n, n>0 \text{ integer}\}$ of *r.v*s generates a max-AR(1) series *if* there exists an innovation sequence $\{\varepsilon_n\}$ of *i.i.d r.v*s such that

$$X_n = \rho X_{n-1} \vee \varepsilon_n, \text{ for some } \rho > 0 \text{ and } \forall \; n>0 \text{ integer}. \tag{1.5}$$

**Theorem.1.1** A sequence $\{X_n\}$ of *r.v*s generates a max-AR(1) series that is marginally stationary (*ie*. $X_n \stackrel{d}{=} X_{n-1} \; \forall \; n>0$ integer) *iff* the distribution of $X_n$ is max-semi-SD($c$), $c=1/\rho$.

**Theorem.1.2** Max-semi-stable($a,b$) laws in (1.1) and (1.2) are max-semi-SD($b$).

With this background here we describe a max-AR(1) model with max-semi-stable marginals and characterize it using semi-selfsimilar extremal processes.

## 2. Results

**Remark.2.1** The max-semi-stable d.fs in $\Phi_{\alpha,a,b}$ can be represented in the form $exp\{-\psi(x)\}$, where $\psi(x)$ satisfies $\psi(x) = a\psi(bx)$, for some $a>1$, $b>1$, and $\alpha>0$ satisfying $ab^{-\alpha} = 1$. Similarly those in $\Psi_{\alpha,a,b}$ can be represented as $exp\{-\psi(x)\}$, where $\psi(x) = a\psi(bx)$, for some $a>1$, $b<1$, and $\alpha>0$ satisfying $ab^{\alpha} = 1$.

**Remark.2.2** In the definition of max-semi-SD laws we do not consider the case $v=0$. Hence (1.3) is equivalent to, $F(x) = F(cx) H(x)\ \forall\ x\in \mathbf{R}$ and for some $c\in (0,1)\cup(1,\infty)$. Thus there is a one-to-one correspondence between the possible scale changes in max-semi-stable and max-semi-SD laws with $v\neq 0$.

**Remark.2.3** Without loss of generality we may consider the range of $b$ as $0<b<1$ in the description of selfsimilarity because (1.4) is equivalent to $\{(b^{-1})^H Y(t)\} \stackrel{d}{=} \{Y(b^{-1}t)\}$ and thus the whole range $b>0$ is covered.

**Theorem.2.1** (A corollary to the theorem.5.1 in Pancheva (2000)). Let $\{Y(t)\}$ be an EP with stationary max-increments. Then $\{Y(t)\}$ is max-semi-stable$(a,b)$ of the extended Frechet type *iff* $\{Y(t)\}$ is $(b^{\alpha}, \frac{1}{\alpha})$-semi-SS. $\{Y(t)\}$ is max-semi-stable$(a,b)$ of the extended Weibull type *iff* $\{Y(t)\}$ is $(b^{-\alpha}, \frac{-1}{\alpha})$-semi-SS.

*Proof.* Let $Y(t)$ be max-semi-stable$(a,b)$ of the extended Frechet type. Since $ab^{-\alpha}=1$,

$$\tfrac{1}{b}Y(at) \stackrel{d}{=} Y(t) \text{ or } Y(at) \stackrel{d}{=} bY(t) \stackrel{d}{=} a^{\frac{1}{\alpha}} Y(t).$$

Hence $Y(at) \stackrel{d}{=} a^{\frac{1}{\alpha}} Y(t)$. Thus $Y(t)$ is $(a, \frac{1}{\alpha})$-semi-SS or $Y(t)$ is $(b^{\alpha}, \frac{1}{\alpha})$-semi-SS.

If $Y(t)$ is max-semi-stable$(a,b)$ of the extended Weibull type, then since $ab^{\alpha}=1$,

$$\tfrac{1}{b}Y(at) \stackrel{d}{=} Y(t) \text{ or } Y(at) \stackrel{d}{=} bY(t) \stackrel{d}{=} a^{\frac{-1}{\alpha}} Y(t).$$

Hence $Y(at) \stackrel{d}{=} a^{\frac{-1}{\alpha}} Y(t)$. Thus $Y(t)$ is $(a, \frac{-1}{\alpha})$-semi-SS or $Y(t)$ is $(b^{-\alpha}, \frac{-1}{\alpha})$-semi-SS.

Converse of both the statements follows from theorem.5.1 in Pancheva (2000).

**Theorem.2.2a** Let $\{Y(t), t \geq 0\}$ be a stationary EP, $X_0 \stackrel{d}{=} Y(1)$ and $\varepsilon_n \stackrel{d}{=} \rho Y(b^\alpha\text{-}1)$, $\forall n$ in (1.5) $b = 1/\rho$. Then (1.5) is marginally stationary with extended Frechet type max-semi-stable$(a,b)$ marginals *if* $\{Y(t)\}$ is $(a, \frac{1}{\alpha})$-semi-SS. Conversely, $\{Y(t)\}$ is $(a, \frac{1}{\alpha})$-semi-SS and the marginals are extended Frechet type max-semi-stable *if* (1.5) is marginally stationary.

*Proof.* Notice that if the *d.f* of $Y(1)$ is $G(u)$ then that of $\rho Y(b^\alpha\text{-}1)$ is $\{G(bu)\}^{b^\alpha - 1}$, $b = 1/\rho$. Further if $\{Y(t)\}$ is $(b^\alpha, \frac{1}{\alpha})$-semi-SS then $Y(b^\alpha t) \stackrel{d}{=} bY(t)$.

Hence $\{G(\frac{u}{b})\} = \{G(u)\}^{b^\alpha}$ and so $\{G(bu)\}^{b^\alpha} = G(u)$, implying $G$ is extended Frechet type. Now under the given assumptions at $n=1$,

$$G_1(u) = G(bu) \{G(bu)\}^{b^\alpha - 1} = \{G(bu)\}^{b^\alpha} = G(u).$$

Thus on iteration (1.5) is marginally stationary with extended Frechet type marginals.

Conversely, let (1.5) is marginally stationary. Then at $n=1$, with $b = 1/\rho$,

$$G(u) = G(bu) \{G(bu)\}^{b^\alpha - 1} = \{G(bu)\}^{b^\alpha}.$$

Hence the marginals and $Y(1)$ are extended Frechet type max-semi-stable$(a,b)$ and consequently the $\{Y(t)\}$ is $(b^\alpha, \frac{1}{\alpha})$-semi-SS. Thus the proof is complete.

**Theorem.2.2b** Let $\{Y(t), t \geq 0\}$ be a stationary EP, $X_0 \stackrel{d}{=} Y(1)$ and $\varepsilon_n \stackrel{d}{=} \rho Y(b^{-\alpha}\text{-}1)$, $\forall n$ in (1.5) $b = 1/\rho$. Then (1.5) is marginally stationary with extended Weibull type max-semi-stable$(a,b)$ marginals *if* $\{Y(t)\}$ is $(a, \frac{-1}{\alpha})$-semi-SS. Conversely, $\{Y(t)\}$ is $(a, \frac{-1}{\alpha})$-semi-SS and the marginals are extended Weibull type max-semi-stable *if* (1.5) is marginally stationary.

*Proof*: Follows as that of Theorem.2.2a.

**Remark.2.4** From theorem.4.2 in Satheesh and Sandhya (2006) we have; Max-stable laws corresponding to (1.1) and (1.2) are max-SD. Consequently analogous to theorems 2.1, 2.2 *a* and *b* we have;

**Theorem.2.3** Let $\{Y(t)\}$ be an EP with stationary max-increments. Then $\{Y(t)\}$ is max-stable of the Frechet type *iff* $\{Y(t)\}$ is $\frac{1}{\alpha}$-SS. $\{Y(t)\}$ is max-stable of the Weibull type *iff* $\{Y(t)\}$ is $\frac{-1}{\alpha}$-SS.

**Theorem.2.4***a* Let $\{Y(t), t \geq 0\}$ be a stationary EP, $X_0 \stackrel{d}{=} Y(1)$ and $\varepsilon_n \stackrel{d}{=} \rho Y(b^{\alpha}-1)$, $\forall n$ in (1.5) $b=1/\rho$. Then (1.5) is marginally stationary with Frechet type max-stable marginals *if* $\{Y(t)\}$ is $\frac{1}{\alpha}$-SS. Conversely, $\{Y(t)\}$ is $\frac{1}{\alpha}$-SS and the marginals are Frechet type max-stable *if* (1.5) is marginally stationary.

**Theorem.2.4***b* Let $\{Y(t), t \geq 0\}$ be a stationary EP, $X_0 \stackrel{d}{=} Y(1)$ and $\varepsilon_n \stackrel{d}{=} \rho Y(b^{-\alpha}-1)$, $\forall n$ in (1.5) $b=1/\rho$. Then (1.5) is marginally stationary with Weibull type max-stable marginals *if* $\{Y(t)\}$ is $\frac{-1}{\alpha}$-SS. Conversely, $\{Y(t)\}$ is $\frac{-1}{\alpha}$-SS and the marginals are Weibull type max-stable *if* (1.5) is marginally stationary.

It may be noted that since in theorems 2.4 *a* and *b* the EPs are SS the statements hold for any $b>0$. Consequently the max-AR(1) series (1.5) holds for all $\rho>0$, including the explosive case of $\rho>1$ as well.